\documentclass[11pt]{amsart}

\usepackage{amsmath,amssymb,latexsym,amsthm,newlfont,enumerate, multirow, hyperref}

\usepackage{graphicx}
\usepackage[all]{xy}

\makeindex
\theoremstyle{plain}

\newtheorem{thm}{Theorem}[section]
\newtheorem{pro}[thm]{Proposition}
\newtheorem{cor}[thm]{Corollary}
\newtheorem{con}[thm]{Conjecture}
\theoremstyle{definition}
\newtheorem{dfn}[thm]{Definition}
\newtheorem{exa}[thm]{Example}
\newtheorem{rem}[thm]{Remark}
\renewcommand{\tilde}{\widetilde}
\DeclareMathOperator{\Sing}{Sing}
\DeclareMathOperator{\Diff}{Diff}
\DeclareMathOperator{\mult}{mult}
\begin{document}
\title{Some examples of Calabi-Yau pairs with maximal intersection and no toric model }
\author{Anne-Sophie Kaloghiros}
\address{Department of Mathematics, Brunel University London, Uxbridge
UB8 3PH, UK}
 \email{anne-sophie.kaloghiros@brunel.ac.uk}
\maketitle

\begin{abstract}It is known that a maximal intersection log canonical Calabi-Yau surface pair is crepant birational to a toric pair. This does not hold in higher dimension: this article presents some examples of maximal intersection Calabi-Yau pairs that admit no toric model. \end{abstract}
\section{Introduction and motivation}

A Calabi--Yau (CY) pair $(X, D_X)$ consists of a normal projective variety $X$ and a reduced sum of integral Weil divisors $D_X$ such that $K_X+D_X\sim_{\mathbb Z} 0$.

The class of CY pairs arises naturally in a number of problems and comprises examples with very different birational geometry. Indeed, on the one hand, a Gorenstein Calabi--Yau variety $X$ can be identified with the CY pair $(X,0)$. On the other hand, if $X$ is a Fano variety, and if $D_X$ is an effective reduced anticanonical divisor, then $(X, D_X)$ is also a CY pair. 

\begin{dfn} 
\begin{enumerate}[(a)]
\item A pair $(X, D_X)$ is (t,dlt) (resp.~(t,lc)) if $X$ is $\mathbb{Q}$-factorial, terminal and $(X, D_X)$ divisorially log terminal (resp.~log canonical). 
\item A birational map $(X, D_X)\stackrel{\varphi} \dashrightarrow (Y, D_Y)$ is volume preserving if $a_E(K_X+D_X)= a_E(K_Y+D_Y)$ for every geometric valuation $E$ with centre on $X$ and on $Y$. 
\end{enumerate}
\end{dfn} 

The dual complex of a dlt pair $(Z, D_Z= \sum D_i)$ is the regular cell complex obtained by attaching an $(|I|-1)$-dimensional cell for every irreducible component of a non-empty intersection $\bigcap_{i\in I}D_i$. 

The dual complex encodes the combinatorics of the lc centres of a dlt pair and \cite{dFKX} show that its PL homeomorphism class is a volume preserving birational invariant.  

By \cite[Theorem 1.9]{CK}, a (t,lc) CY pair $(X, D_X)$ has a volume preserving (t,dlt) modification $(\widetilde X, D_{\tilde X}) \to (X, D_X)$, and the birational map between two such modifications is volume preserving. 

Abusing notation, I call dual complex the following volume preserving birational invariant of a (t,lc) CY pair $(X, D_X)$. 

\begin{dfn} $\mathcal{D}(X, D_X)$ is the PL homeomorphism class of the dual complex of a volume preserving (t,dlt) modification of $(X, D_X)$. 
\end{dfn}

As the underlying varieties of CY pairs range from CY to Fano varieties, they can have very different birational properties. However, $X$ being Fano is not a volume preserving birational invariant of the pair $(X, D_X)$. Following \cite{KX}, I consider the following volume preserving birational invariant notion:
\begin{dfn}
A (t,lc) CY pair $(X, D_X)$ has maximal intersection if $\dim \mathcal D(X, D_X)= \dim X-1$.\end{dfn}

In other words, $(X, D_X)$ has maximal intersection if there is a volume preserving (t,dlt) modification of $(X, D_X)$ with a $0$-dimensional log canonical centre. Maximal intersection CY pairs have some Fano-type properties; Koll\'ar and Xu show the following:
\begin{thm} Let $(X, D_X)$ be a dlt maximal intersection CY pair, then:
\begin{enumerate}[1.]\item \cite[Proposition 19]{KX}
 $X$ is rationally connected,
\item\cite[Theorem 21]{KX}  there is a volume preserving map $(X, D_X)\stackrel{\varphi}\dashrightarrow (Z, D_Z)$ such that $D_Z$ fully supports a big and semiample divisor. 
\end{enumerate}\end{thm}

\begin{rem} 

The expression ``Fano-type" should be understood with a pinch of salt. Having maximal intersection is a degenerate condition: a general (t,lc) CY pair $(X,D_X)$ with $X$ Fano and $D_X$ a reduced anticanonical section needs not have maximal intersection. 
\end{rem}
\begin{dfn} A toric pair $(X, D_X)$ is a (t,lc) CY pair formed by a toric variety and the reduced sum of toric invariant divisors. 

A toric model is a volume preserving birational map to a toric pair. \end{dfn}
\begin{exa}\label{ex1}
A  CY pair with a toric model has maximal intersection. 
\end{exa}
\begin{rem}\label{rem1} In dimension $2$, the converse holds: maximal intersection CY surface pairs are precisely those with a toric model \cite{GHK}. 
\end{rem}

The characterisation of CY pairs with a toric model is an open and difficult problem. A characterisation of toric pairs was conjectured by Shokurov and is proved in \cite{BMcKSZ}, but it is not clear how to refine it to get information on the existence of a toric model.  
A motivation to better understand the birational geometry of CY pairs and their relation to toric pairs comes from mirror symmetry. 

The mirror conjecture extends from a duality between Calabi-Yau varieties to a correspondence between Fano varieties and Landau-Ginzburg models, i.e.~non-compact K\"ahler manifolds endowed with a superpotential. 
Most known constructions of mirror partners rely on toric features such as the existence of a toric model or of a toric degeneration. 
In an exciting development, Gross, Hacking and Keel conjecture the following construction for mirrors of maximal intersection CY pairs.    

\begin{con} \label{msconj}\cite{GHK} Let $(Y,D_Y)$ be a simple normal crossings maximal intersection CY pair. Assume that $D_Y$ supports an ample divisor, let $R$ be the ring $k[\mathrm{Pic}(Y)^\times]$, $\Omega$ the canonical volume form on $U$ and \[ 
U^{\mathrm{trop}}(\mathbb Z)= \Big\{\mbox{divisorial valuations } v \colon k(U)\setminus \{0\} \to \mathbb Z  \mbox{ with } v(\Omega)<0\Big\}\cup \{0\}.\] 
Then, the free $R$-module $V$ with basis $U^{\mathrm{trop}}(\mathbb Z)$ has a natural finitely generated $R$-algebra structure whose structure constants are non-negative integers determined by counts of rational curves on $U$. 

Denote by $K$ the torus $\mathrm{Ker} \{\mathrm{Pic} (Y) \to \mathrm{Pic} (U)\}$. The fibration $$p\colon \mathrm{Spec}(V )\to \mathrm{Spec}(R) = T_{\mathrm{Pic}(Y)}$$ is a $T_K$-equivariant flat family of affine maximal intersection log CY varieties. The quotient $$\mathrm{Spec}(V )/T_K \to T_{\mathrm{Pic}(U)}$$ only depends on $U$ and is the mirror family of $U$.
\end{con}

Versions of Conjecture~\ref{msconj} are proved for cluster varieties in \cite{GHKK}, but relatively few examples are known.

The goal of this note is to present examples of maximal intersection CY pairs that do not admit a toric model and for which one can hope to construct the mirror partner proposed in Conjecture~\ref{msconj} (see Section~\ref{3f} for a precise statement).

\section{Auxiliary results on $3$-fold CY pairs}
\label{3f}

The examples in Section~\ref{exa} are $3$-fold maximal intersection CY pairs whose underlying varieties are birationally rigid. In particular, such pairs admit no toric model; this shows that \cite{GHK}'s results on maximal intersection surface CY pairs do not extend to higher dimensions.  
In this section, I first recall some results on birational rigidity of Fano $3$-folds. Then, I introduce the (t,dlt) modifications suited to the construction outlined in Conjecture~\ref{msconj} and discuss the singularities of the boundary $D_X$.   
\subsection{Birational rigidity} Let $X$ be a terminal $\mathbb Q$-factorial Fano $3$-fold. When $X$ has Picard rank $1$, $X$ is a Mori fibre space, i.e.~an end product of the classical MMP. 

\begin{dfn} 
A birational map $Y/S\stackrel{\varphi}\dashrightarrow Y'/S'$ between Mori fibre spaces $Y/S$ and $Y'/S'$ is square if it fits into a commutative square 
\[
\xymatrix{ Y\ar[d]\ar@{-->}[rr]^\varphi&\quad & Y'\ar[d] \\
S \ar@{-->}[rr]^g& & S'
}
\]
where $g$ is birational and the restriction $Y_{\eta}\stackrel{\varphi_{\eta}} \dashrightarrow Y'_{\eta}$ is biregular, where $\eta$ is the function field of the base $k(S)$. 

A Mori fibre space $Y/S$ is (birationally) rigid if for every birational map $Y/S\stackrel{\varphi}\dashrightarrow Y'/S'$ to another Mori fibre space, there is a birational self map $Y/S\stackrel{\alpha}\dashrightarrow X/S$ such that $\varphi \circ \alpha$ is square.
\end{dfn}

In particular, if $X$ is a rigid Mori fibre space, then $X$ is non-rational and no (t,lc) CY pair $(X,D_X)$ admits a toric model.

Non-singular quartic hypersurfaces $X_4\subset \mathbb P^4$ are probably the most famous examples of birationally rigid $3$-folds \cite{IM}.   
Some mildly singular quartic hypersurfaces are also known to be birationally rigid, in particular, we have:
 \begin{pro}\cite{Che,Me04}\label{rigid}
 Let $X_4\subset \mathbb P^4$ be a quartic hypersurface with no worse than ordinary double points. If $|\Sing(X) |\leq 8$, then $X$ is $\mathbb Q$-factorial (in particular, $X$ is a Mori fibre space) and is birationally rigid. 
 \end{pro}

\subsection{Singularities of the boundary}
\label{sings}
I now state some results on the singularities of the boundary of a $3$-fold (t,lc) CY pair.
Let $(X, D_X)$ be a $3$-fold (t,lc) CY pair and $(\widetilde X, D_{\tilde X})$ a (t,dlt) modification. A stratum of $(\widetilde X, D_{\tilde X})$ is an irreducible component of a non-empty intersection of components of $D_{\tilde X}$.
Given a stratum $W$, there is a divisor $\Diff_WD_{\tilde X}$ on $W$ such that 
$(W, \Diff_W D_{\tilde X})$ is a lc CY pair and $$K_W+ \Diff_W D_{\tilde X}\sim_{\mathbb Q} \big(K_{\tilde X}+ D_{\tilde X}\big)_{|W}.$$

When $K_{\tilde X}+ D_{\tilde X}$ is Cartier and $D_{\tilde X}$ reduced, $\Diff_W D_{\tilde X}$ is the sum of the restrictions of the components of $D_{\tilde X}$ that do not contain $W$. 

In particular, for any irreducible component $S$ of $D_{\tilde X}$, the link of $[S]$ in $\mathcal D(X,D_X)$ is the dual complex $\mathcal D(S, \Diff_S D_{\tilde X})$. Therefore, if 
$(X, D_X)$ has maximal intersection, so does $(S, \Diff_SD_{\tilde X})$. By the results of \cite{GHK}, $(S, \Diff_SD_X)$ then has a toric model. 

As $X$ has terminal singularities, $X$ is normal and Cohen-Macaulay. Any Cartier component $S$ of the boundary $D_X$ is Cohen-Macaulay and satisfies Serre's condition $S_2$. By \cite[Proposition 16.9]{FA92}, $(S, \Diff_SD_X)$ is semi log canonical (slc). In particular, if $X$ is Gorenstein and $D_X$ irreducible, $D_X$ has slc singularities. 

I am particularly interested in producing examples of (t,lc) CY pairs for which the mirror partners proposed in Conjecture~\ref{msconj} (see also \cite{GS}) can be constructed; this motivates the following definition:
\begin{dfn}\label{good} 
A (t,dlt) modification $(\widetilde X, D_{\tilde X})\to (X,D_X)$ is called good if $(\widetilde X, D_{\tilde X})$ is log smooth in the sense of log geometry, that is if the components of $D_{\tilde X}$ are non-singular and if $\widetilde X$ has only cyclic quotient singularities. \end{dfn}

An immediate consequence of the definition is that if $(\widetilde X, D_{\tilde{X}})\stackrel{f}\to (X, D_X)$ is a good (t,dlt) modification and $D_X= \sum_{i}D_i$, then \[ D_{\tilde X}= \sum_if_*^{-1}D_i+E,\]where $E$ is reduced and $f$-exceptional, and the restriction of $f$ to $f_*^{-1}D_i$ is a resolution for all $i$.

{\bf Normal singularities}
Let $p\in \Sing(D_i)$ be an isolated singularity lying on a single component of the boundary. The restriction $f_i\colon \widetilde D_i\to D_i$ is a resolution and we have: 
 \[ K_{\tilde D_i}= (K_{\tilde X}+ D_i)_{\vert \tilde D_i} = (f_{|\tilde D_i})^*K_{D_i}-(E)_{|\tilde D_i}\] 
 where $E$ is defined by $K_{\tilde X}+ f_*^{-1}D_{\tilde X}+E= f^*(K_X+D_X)$.  
 
We now assume that $D_i$ is Cartier, as is the case when $X$ is Gorenstein and $D_X$ irreducible. Without loss of generality, assume that $\Sing(D_i)= p$. Then, $p$ is canonical if $E\cap \widetilde D_i= \emptyset$, and elliptic otherwise. Indeed, let $$f_i \colon \widetilde D_i\stackrel{q}\to \overline D_i\stackrel{\mu}\to D_i$$ be the factorisation through the minimal resolution of $(p\in D_i)$. Then, $q$ is either an isomorphism or an isomorphism at the generic point of each component of $E_{|\tilde D_i}$ because $f$ is volume preserving. We have: $K_{\overline D_i}= \mu^*K_{D_i}-Z$, where the effective cycle $Z= q_*(E_{D_i})$ is either empty (and $p$ is canonical) or a reduced sum of $\mu$-exceptional curves (and $p$ is elliptic). In the second case, $Z\sim -K_{\widetilde D_i}$ is the fundamental cycle of $(p\in D_i)$. If $Z$ is irreducible, it is reduced and has genus $1$; if not, every irreducible component of $Z$ is a smooth rational curve of self-intersection $-2$. 

When $p$ is elliptic,  $Z$ is reduced and $p$ is a Kodaira singularity \cite[Theorem 2.9]{Kar}, i.e.~a resolution is obtained by blowing up points of the singular fibre in a degeneration of elliptic curves; further, in Arnold's terminology, the singularity $p$ is uni or bimodal. 

Further, $p\in D_i$ is a hypersurface singularity (resp.~a codimension $2$ complete intersection, resp.~not a complete intersection) when $-3\leq Z^2\leq -1$ (resp.~$Z^2= -4$, resp.~$Z^2\leq -5$) \cite{Lau}. When $-1\leq Z^2\leq -4$, normal forms are known for $p\in D_i$: Table 1 lists normal forms of slc hypersurface singularities, while normal forms of codimension $2$ complete intersections elliptic singularities are given in \cite{Wall}. \begin{table}[htbp]\label{table}
\label{slchyp}
\centering
\small
\begin{tabular}{llllcc}
\hline
type & name & symbol & equation $f\in \mathbb C[x,y,z]$ & &$\mult_0f$\\
\hline
\noalign{\smallskip}
terminal & smooth & $A_0$ & $x$ && 1\\
\noalign{\smallskip}
\hline 
\noalign{\smallskip}
\multirow{5}{*}{canonical} & \multirow{5}{*} {du Val} &
$A_n$ & $x^2+y^2+ z^{n+1}$  &$n\geq 1$ & 2 \\
&&$D_n$ & $x^2+z(y^2+ z^{n-2})$ &$n\geq 4$ & 2\\
&&$E_6$ & $x^2+y^3+ z^4$ && 2\\
&&$E_7$ & $x^2+y^3+ yz^3$ && 2\\
&&$E_8$ & $x^2+y^3+ z^5$ && 2\\
\noalign{\smallskip}
\hline
\noalign{\smallskip}
\multirow{ 4}{*}{lc} &  
\multirow{3}{*}{simple elliptic} 
& $X_{1,0}$ & $x^2+ y^4+z^4+ \lambda xyz$& $\lambda^4\neq 64$& 2\\
&& $J_{2,0}$ & $x^2+ y^3+z^6+ \lambda xyz$ &$\lambda^6\neq 432$& 2\\

&&$T_{3,3,3}$ & $x^3+ y^3+z^3+ \lambda xyz$& $\quad \lambda^3\neq -27$& 3\\ \noalign{\smallskip}\cline{2-6}
\noalign{\smallskip}

&cusp & $T_{p,q,r}$ & $x^p+ y^q+z^r+ xyz$ &$\frac1{p}+\frac1{q}+\frac1{r}<1$& 2 or 3\\ \noalign{\smallskip}\hline\noalign{\smallskip}

\multirow{7}{*}{slc} & normal crossing & $A_{\infty}$ & $x^2+y^2$ &&2\\ \cline{2-6} \noalign{\smallskip}
& pinch point & $D_{\infty}$ & $x^2+y^2z$ &&2\\ \noalign{\smallskip}\cline{2-6} \noalign{\smallskip}

&\multirow{5}{*}{degenerate cusp} & $T_{2, \infty, \infty}$ & $x^2+y^2+z^2$ && 2\\
&& $T_{2, q, \infty}$ & $x^2+y^2(z^2+ y^{q-2})$ &$q\geq 3$& 2\\
&& $T_{\infty, \infty, \infty}$ & $xyz$ && 3\\
&& $T_{p,  \infty, \infty}$ & $xyz+ x^p$ &$p\geq 3$& 3\\
&& $T_{p, q, \infty}$ & $xyz+ x^p+ y^q$ &$q\geq p\geq 3$& 3\\ \noalign{\smallskip} \hline
\end{tabular}
\smallskip
\caption{Dimension $2$ slc hypersurface singularities}
\end{table}

\section{Examples of rigid maximal intersection $3$-fold CY pairs}\label{exa}

All the examples below are $(t,lc)$ CY pairs $(X, D_X)$ which admit no toric model. Except for Example~\ref{exa5}, all underlying varieties $X$ are birationally rigid quartic hypersurfaces by Proposition~\ref{rigid}; the underlying variety in Example~\ref{exa5} is a smooth cubic $3$-fold, and therefore non-rational.

\subsection{Examples with normal boundary}
\
\begin{exa}\label{exa2} Consider the CY pair $(X, D_X)$ where $X$ is the nonsingular quartic hypersurface
\[ X= \{ x_1^4+ x_2^4+x_3^4+ x_0x_1x_2x_3+ x_4(x_0^3+x_4^3)=0\}\] and $D_X$ is its hyperplane section $X\cap \{ x_4=0\}$.

The quartic surface $D_X$ has a unique singular point $p=(1{:}0{:}0{:}0{:}0)$, and using the notation of Table 1, $p$ is locally analytically equivalent to a  $T_{4,4,4}$ cusp $$\{0\}\in \{ x^4+y^4+z^4 + xyz=0\}.$$ $D_X$ is easily seen to be rational: the projection from the triple point $p$ is \[D_{X}\dashrightarrow \mathbb P^2_{x_1,x_2, x_3};\] this map is  the blowup of the $12$ points $\{ x_1^4+x_2^4+x_3^4=x_1x_2x_3=0\}$, of which $4$ lie on each coordinate line $L_i= \{ x_i=0\}$, for $i=1,2,3$. 

I treat this example in detail and construct explicitly a good (t,dlt) modification of the pair $(X, D_X)$. 

Let $f\colon X_p\to X$ be the blowup of $p$, then $X_p$ is non-singular, the exceptional divisor $E$ satisfies $(E, \mathcal O_E(E))= (\mathbb P^2,\mathcal O_{\mathbb P^2}(-1))$, and if $D$ denotes the proper transform of $D_{X}$, we have:
\[ K_{X_p}+ D+E = f^*(K_X+D).\]

Explicitly, the blowup $\mathcal F\to \mathbb P^4$ of $\mathbb P^4$ at $p$ is the rank $2$ toric variety $\mathrm{TV}(I, A)$, where $I=(u,x_0)\cap (x_1, \dots , x_4)$ is the irrelevant ideal of $\mathbb C[u, x_0, \dots, x_4]$ and $A$ is the action of $\mathbb C^\ast \times \mathbb C^\ast$ with weights: 
\begin{equation}
\left( \begin {array}{cccccccc}
 u & x_0& s_1&s_2 & s_3&s_4\\
1 & 0 & -1 & -1 &-1 &-1 \\
 0 & 1 & 1& 1& 1& 1\end{array}
\right).
\end{equation}

The equation of $X_p$ is
\[ X_p=\{ u^2\big(u(s_1^4+s_2^4+s_3^4+ (x_0s_1s_2s_3)\big)+ s_4(x_0^3+u^3s_4^3)=0\}, \]
while $E= \{ u=0\}$ and $D= \{ u(s_1^4+s_2^4+s_3^4)+ x_0s_1s_2s_3=0\}$. By construction, $E$ is the projective plane with coordinates $s_1, s_2, s_3$. Note that $(X_p, D+E)$ is not dlt because $D\cap E= \{ x_0s_1s_2s_3=0\}$ consists of $3$ concurrent lines $C_1, C_2, C_3$.

Consider $g_1\colon X_1\to X_p$ the blowup of the nonsingular curve $$C_1= \{ u=s_1=s_4=0\}\subset X_p.$$ The exceptional divisor of $g_1$ is a surface $E_1\simeq \mathbb P(\mathcal N_{C_1/ X_p})$, and since $C_1\simeq \mathbb P^1$, the restriction sequence of normal bundles gives 
\[ \mathcal N_{C_1/ X_p} \simeq \mathcal N_{C_1/E}\oplus (\mathcal N_{E/X_p})_{|C_1}\simeq \mathcal O_{C_1}(1)\oplus {\mathcal O_{E}(-1)}_{|C_1},\] so that $E_1= \mathbb F_2$. Further, 
\[ K_{X_1}+ D+E+E_1= g_1^*(K_{X_{p}}+D+E)\]
where, abusing notation, I denote by $D$ and $E$ the proper transforms of the divisors $D$ and $E$. The ``restricted pair" on $E_1$ is a surface CY pair $(E_1, (D+E)_{|E_1})$ by adjunction. By construction, $E\cap E_1$ is the negative section $\sigma$. The curve $\Gamma= D\cap E_1$ is irreducible, and since $(D+E)_{|E_1}$ is anticanonical, we have   
\[ \Gamma\sim \sigma + 4f \mbox{ where $f$ is a fibre of } \mathbb F_2\to \mathbb P^1, \mbox{ and } \Gamma^2= 6, \Gamma\cdot E_{|E_1}=2.\]
The divisors $D,E,E_1$ meet in two points, the dual complex $\mathcal D(X_1, D+E+E_1)$ is not simplicial it is a sphere $S^2$ whose triangulation is given by $3$ vertices on an equator. While not strictly necessary, we consider a further blowup to obtain a (t,dlt) pair with simplicial dual complex. 

Denote by $C_2$ the proper transform of the curve \[ \{ u=s_2=s_4=0\}.\] Then $C_2\subset E\cap D$ is rational, and as above $$\mathcal N_{C_2/X_1}\simeq \mathcal N_{C_2/E}\oplus({\mathcal N}_{E/X_2})_{|C_2}= \mathcal O_{C_2}(1)\oplus \mathcal O_{C_2}(-2).$$ Let $g_2\colon X_2\to X_1$ be the blowup of $C_2$, then the exceptional divisor of $g_2$ is a Hirzebruch surface \[ E_2\simeq \mathbb P_{\mathbb P^1}(\mathcal N_{C_2/X_1})\simeq\mathbb F_3.\] Still denoting by $D,E,E_1$ the strict transforms of $D,E,E_1$, we have:
\[
K_{X_2}+ D+E+E_1+E_2= g_2^*(K_{X_1}+ D+E+E_1).\]
The pair $(X_2, D+E+E_1+E_2)$ is dlt; the composition \[ g_2\circ g_1\circ f\colon (\widetilde X, D_{\tilde X})= (X_2, D+E+E_1+E_2)\to (X, D_{X})\] is a good (t,dlt) modification. 

The ``restrictions" of $(\widetilde X, D_{\tilde X})$ to the component of the boundary are the following surface anticanonical pairs:
\begin{enumerate}[-]
\item On $D$: $ (E+E_1+E_2)_{|D}$ is a cycle of $(-3)$-curves, the morphism $D\to D_{X}$ is the familiar resolution of the $T_{4,4,4}$ cusp singularity;
\item On $E$: $(D+E_1+E_2)_{E}$ is the triangle of coordinate lines with self-intersections $(1,1,1)$;
\item On $E_1$: $(D+E+E_2)_{E_1}$ is an anticanonical cycle with self-intersections $(5,-3,-1)$;
\item On $E_2$:  $(D+E+E_1)_{E_2}$ is an anticanonical cycle with self-intersections $(5,-3,0)$ (as above, $E_{|E_2}\sim \sigma$ is a negative section, ${E_1}_{|E_2}\sim f$ a fibre of $\mathbb F_3\to \mathbb P^1$, and $D_{|E_2}\sim 4f+\sigma$).
\end{enumerate}  
It follows that the dual complex $\mathcal D(X, D_X)$ is PL homeomorphic to a tetrahedron and $(X, D_X)$ has maximal intersection. Note that $(0\in D_X)$ is a maximal intersection lc point, and since $D_X$ is a rational surface, it has a toric model. 
\end{exa}

\begin{exa}
Let $X$ be the hypersurface
\[ X= \{ x_3(x_0^3+ x_1^3)+ x_2^4+ x_0x_1x_2x_3+ x_4(x_3^3+x_4^3)=0\}, \] and $D_X$ its hyperplane section $X\cap \{ x_4=0\}$.

The quartic $X$ has $3$ ordinary double points at the intersection points $$L \cap \{x_0^3+ x_1^3=0\},$$ where $L$ is the line $\{ x_2=x_3=x_4=0\}$. 
The singular locus of $D_X$ is $\Sing(X)\cup \{p\}$, where $p =(0{:}0{:}0{:}1{:}0)$ is a $T_{3, 3, 4}$ cusp, i.e.~locally analytically equivalent to \[ \{0\}\in \{ x^3+y^3+z^4+ xyz=0\}.\]  
The quartic surface $D_X$ is rational; the projection of $D_X$ from $p$ is \[D_{X}\dashrightarrow \mathbb P^2_{x_0,x_1, x_2};\] this map is defined outside of the $12$ points (counted with multiplicity) defined by $\{ x_2^4= x_0^3+x_1^3+ x_0x_1x_2=0\}$. 

If $\tilde X\stackrel{f}\to X$ is the composition of the blowups at the ordinary double points and at $p$, $\widetilde X$ is smooth and $D_{\tilde X}$ is non-singular, so that $f$ is a good (t,dlt) modification. 

The minimal resolution of $p\in D_X$ is a rational curve with self intersection $C^2=-3$. Explicitly, taking the blowup of $X$ at $p$, the proper transform is a rational surface $D$. The exceptional curve is the preimage of a nodal cubic in $\mathbb P^2$ blown up at 12 points counted with multiplicities. Note that $(\widetilde X, D+E)$ is not dlt, but in order to obtain a (t,dlt) modification, we just need to blowup the node of $D\cap E$ which is a nonsingular point of $\widetilde X, D$ and $E$. The (t,dlt) modification of $(X,D_X)$ in a neighbourhood of $p$ is good and the associated dual complex is $2$-dimensional. 

The pair $(X, D_X)$ has maximal intersection; but as in the previous examples, $X$ is rigid, so that $(X, D_X)$ can have no toric model.  
\end{exa}
\begin{exa} 
Let $X$ be the nonsingular quartic hypersurface
\[ X= \{ x_0^3x_3+x_1^4+x_2^4+ x_0x_1x_2x_3+ x_4(x_3^3+x_4^3)=0\}\subset \mathbb P^4\]
and $D_X$ its hyperplane section $X\cap \{ x_4=0\}$.

The surface $D_X$ has a unique singular point $p=(0{:}0{:}0{:}1{:}0)$ of $D_X$, which is a cusp $T_{3,4,4}$, i.e.~is locally analytically equivalent to \[\{0\}\in \{ x^3+y^4+z^4 + xyz=0\}.\] 

As in Example~\ref{exa2}, $X$ is non-singular, and finding a good (t,dlt) modification of $(X, D_X)$ will amount to taking a minimal resolution of the singular point of $D_X$. Let $X_p\to X$ be the blowup of $X$ at $p$; $X_p$ is non-singular and if $D$ denotes the proper transform of $D_X$, and $E$ the exceptional divisor, $D\cap E$ consists of $2$ rational curves of self intersection $-3$ and $-4$. These curves are the proper transforms of $\{x_0=0\}$ and of $\{ x_0^2+ x_1x_2\}$ under the blow up of $\mathbb P^2_{x_0,x_1,x_2}$ at the points 
\[\{ x_1^4+x_2^4= x_0(x_1x_2+ x_0^2)=0\}.\]

 The dual complex consists of $3$ vertices that are joined by edges and span $2$ distinct faces: $\mathcal D(X, D_X)$ is PL homeomorphic to a sphere $S^2$ whose triangulation is given by $3$ vertices on an equator. The CY pair $(X, D_X)$ has maximal intersection but no toric model.
\end{exa}
\subsection{Examples with non-normal boundary}

\begin{exa}\label{exa5} This example is due to R.~Svaldi. 
Consider the cubic $3$-fold 
\[ X= \{ x_0x_1x_2+ x_1^3+x_2^3+ x_3q+x_4q'=0\}\subset \mathbb P^4\]
where $q, q'$ are homogeneous polynomials of degree $2$ in $x_0, \cdots, x_4$. If the quadrics $q$ and $q'$ are general and if \[(q(1,0,0,0,0), q'(1,0,0,0,0))\neq (0,0),\] then $X$ and $S= \{ x_3=0\}\cap X$ and $T= \{ x_4=0\}\cap X$ are nonsingular. 

Let $D_X$ be the anticanonical divisor $S+T$. 
The curve $C=S\cap T= \Pi\cap X$ for $\Pi= \{ x_3=x_4=0\}$ is a nodal cubic. It follows that both $(S, C)$ and $(T, C)$ are log canonical, and therefore so is $(X, D_X)$. 

Since $S$ and $T$ are smooth, $\Sing(D_X)= S\cap T= C$, and if $p$ is the node of $C$, we have:
\begin{eqnarray*}
(p\in D_X) \sim \{0\}\in \big\{ (xy+x^3+y^3+z)(xy+x^3+y^3+t)=0\big\}\\
\sim \{0\}\in \big\{(xy+z)(xy+t)=0\big\}\sim \{0\}\in\big\{ (xy+z)(xy-z)=0\big\}.\end{eqnarray*}

Thus, $p\in D_X$ is a double pinch point, i.e.~$p$ is locally analytically equivalent to $\{0\}\in \{ x^2y^2-z^2=0\}$. 

We now construct a good (t,dlt) modification of $(X, D_X)$.  
Let $f\colon X_C\to X$ be the blowup of $X$ along $C$; $\Sing(X_C)$ is an ordinary double point. 

Indeed, let $\Pi= \{ x_3=x_4=0\}$, then $f$ is the restriction to $X$ of the blowup  $\mathcal F\to \mathbb P^4$, where $\mathcal F$ is the rank $2$ toric variety $\mathrm{TV}(I, A)$, where $I=(u,x_0, x_1, x_2)\cap (x_3, x_4)$ is the irrelevant ideal of $\mathbb C[u, x_0, \dots, x_4]$ and $A$ is the action of $\mathbb C^\ast \times \mathbb C^\ast$ with weights: 
\begin{equation*}
\left( \begin {array}{cccccc}
 u & x_0& x_1&x_2 & x_3&x_4\\
1 & 0 & 0 & 0 & -1 &-1 \\
 0 & 1 & 1& 1& 1& 1\end{array}
\right),
\end{equation*}
The equation of $X_C$ is $$\{ x_0x_1x_2+ x_1^3+x_2^3+ u(x_3q+x_4q')=0\},$$ so that $X_C$ has a unique singular point at $$x_0-1=u= x_1=x_2= x_3q(1,0,0,0,0)+x_4q'(1,0,0,0,0)=0,$$ and this is a $3$-fold ordinary double point. 
In addition, denoting by $E_f= \{ u=0\}\cap X_C$ the exceptional divisor, we have
\[ K_{X_C}+ \tilde S+ \tilde T+ E_f= K_X+S+T,\]
so that the pair $(X_C, \tilde S+ \tilde T+E_f)$ is a (t,lc) CY pair. 

The pair $(X_C, \tilde S+ \tilde T+E_f)$ is not dlt as the boundary has multiplicity $3$ along the fibre $F$ over the node of $S\cap T$. The blowup of $F$ is not $\mathbb Q$-factorial, therefore in order to obtain a good (t,dlt) modification, we consider the divisorial contraction $g\colon \tilde X\to X_C$ centred along $F$. This is obtained by (a) blowing up the node, (b) then blowing up the proper transform of $F$, (c) flopping a pair of lines with normal bundle $(-1,-1)$ and (d) contracting the proper transform of the $\mathbb P^1\times \mathbb P^1$ above the node to a point $\frac1{2}(1,1,1)$. The exceptional divisor of $g$ is denoted by $E_g$. 

The pair $(\tilde X, \tilde S+ \tilde T+ \tilde E_f+E_g)$ is the desired (t,dlt) modification of $(X, D_X)$, and it has maximal intersection. The dual complex is PL homeomorphic to a tetrahedron.  

\end{exa}

\begin{exa}
Let $X$ be the quartic hypersurface
\[ X= \{ x_1^2x_2^2+ x_1x_2x_3l+ x_3^2q+ x_4f_3=0\}\subset \mathbb P^4,\] 
where $l$ (resp.~$q$) is a general linear (resp.~quadratic) form in $x_0, \cdots, x_3$, and $f_3$ a general homogeneous form of degree $3$ in $x_0, \cdots, x_4$. Let $D_X$ be the hyperplane section $X\cap \{ x_4=0\}.$

As $l, q$ and $f_3$ are general, $X$ has $6$ ordinary double points. Indeed, denote by $L=\{ x_1=x_3=x=4=0\}$ and $L'=\{ x_2=x_3=x=4=0\}$, then 
\[ \Sing(X) =\big\{ L\cap\{f_3=0\}\big\} \cup \big\{L'\cap \{f_3=0\}\big\}= \{q_1, q_2, q_3\}\cup  \{ q'_1, q'_2, q'_3\}\]
which consists of $3$ points on each of the lines. 
In the neighbourhood of each point $q_i$ (resp.~$q'_i$) for $i=1,2,3$, the equation of $X$ is of the form 
\[ \{0\}\in \{ xy+zt=0\}\]
(and $D_X= \{t=0\}$) so that all singular points of $X$ are ordinary double points. 
The quartic hypersurface $X$ is birationally rigid by Proposition~\ref{rigid}. 

The surface $D_X$ is non-normal as it has multiplicity $2$ along $L$ and $L'$. The point $p= L\cap L'$ is locally analytically equivalent to $$\{0\}\in \{x^2y^2+z^2=0\},$$ so that $p\in D_X$ is a double pinch point. We conclude that the surface $D_X$ has slc singularities, and hence $(X, D_X)$ is a (t,lc) CY pair. 

We construct a good (t,dlt) modification as follows. 

First, since $\Sing(X) \cap L$ (resp.~$\Sing (X) \cap L'$) is non-empty, the blowup of $X$ along $L$ (resp.~along $L'$) is not $\mathbb Q$-factorial. In order to remain in the (t,dlt) category, we consider the divisorial extraction $f\colon X_L\to X$ centered on $L$ (resp.~$L'$).This is obtained by (a) blowing up the $3$ nodes lying on $L$, (b) blowing up the proper transform of $L$, (c) flopping $3$ pairs of lines with normal bundle $(-1, -1)$ and (d) contracting the proper transforms of the three exceptional divisors $\mathbb P^1\times \mathbb P^1$ lying above the nodes to points $\frac1{2}(1,1,1)$. The exceptional divisor of $f$ is denoted by $E$. 
Let $p\colon \widetilde X \to X$ denote the morphism obtained by composing the divisorial extraction centered on $L$ with that centered on $L'$ (in any order), and let $E,E'$ denote the exceptional divisors of the divisorial extractions. Then
\[ K_{\widetilde X}+ \widetilde D+ E+E'= p^*(K_X+D)\]
is a (t,dlt) modification of $(X, D_X)$ and it has maximal intersection. The dual complex $\mathcal D(X, D_X)$ is  is PL homeomorphic to a sphere $S^2$ whose triangulation is given by $3$ vertices on an equator. \end{exa}

\section{Further results on quartic $3$-fold CY pairs: beyond maximal intersection} \label{4ics}

This section concentrates on (t,lc) CY pairs $(X,D_X)$, where $X$ is a factorial quartic hypersurface in $\mathbb P^4$ and $D$ is an irreducible hyperplane section of $X$. I give some more detail on the possible dual complexes of such pairs. 
 
As explained in Section~\ref{sings}, $D_X$ is slc because $(X,D_X)$ is lc. 
In order to study completely the dual complexes of such (t,lc) CY pairs, one needs a good understanding of the normal forms of slc singularities that can lie on $D$. In the case of a general Fano $X$, this step would require additional work, but here, $D_X$ is a quartic surface in $\mathbb P^3$ and the study of singularities of such surfaces has a rich history. I recall some results directly relevant to the construction of degenerate $CY$ pairs $(X,D_X)$. 
The classification of singular quartic surfaces in $\mathbb P^3$ can be broken in three independent cases. 
\begin{enumerate}[(a)]
\item Quartic surfaces with no worse than rational double points: the minimal resolution is a $K3$ surface. Possible configurations of canonical singularities were studied by several authors using the moduli theory of $K3$ surfaces;  there are several thousands possible configurations. The pair $(X,D_X)$ is (t,dlt) and the dual complex of $(X,D_X)$ is reduced to a point. 
\item Non-normal quartic surfaces were classified by Urabe \cite{Ur}; there are a handful of cases recalled in Theorem~\ref{nonnormal4ics}.    
\item Non-canonical quartic surfaces with isolated singularities. These are studied by Wall \cite{Wallp} and Degtyarev \cite{Deg} among others; their results are recalled in Theorem~\ref{isosingquartic}.  
\end{enumerate}

\begin{thm}\label{nonnormal4ics}\cite{Ur}A non-normal quartic surface $D\subset \mathbb P^3$ is one of:
\begin{enumerate}
\item[1.]
the cone over an irreducible plane quartic curve with a singular point of type $A_1$ or $A_2$.
\item[2.] a ruled surface over a smooth elliptic curve $G$, $D= \varphi_{\mathcal L}(Z)$, where: 
\item[(a)] $\mathcal L= \mathcal O_Z(C_1)\otimes \pi^*M$, and $Z= \mathbb P_G(\mathcal O_G\oplus N)$, for \subitem -$M$ a line bundle of degree $2$ and \subitem -$N$ a non-trivial line bundle of degree $0$.

Denoting by $L_i$ the images by $\varphi_{\mathcal L}$ of the sections of $Z$ associated to $\mathcal O_G\oplus N \to \mathcal O_G$ and $\mathcal O_G\oplus N \to N$, 
$\Sing(D)= L_1\cup L_2$. 
 \item[(b)] $\mathcal L= \mathcal O_Z(C)\otimes \pi^*M$, and $Z= \mathbb P_G(E)$ for a rank $2$ vector bundle $E$ that fits in a non-splitting  
 \[0 \to \mathcal O_G\to E \to \mathcal O_G \to 0.\]
 Denoting by $L$ the image by $\varphi_{\mathcal L}$ of a section $G \to Z$, $\Sing D=L$.
\item[3.] a rational surface $D\subset \mathbb P^3$ which is
\item[(a)] the image of a smooth $S\subset \mathbb P^5$ under the projection from a line disjoint from $S$; $D$ has no isolated singular point and 
\subitem-$S= v_2(\mathbb P^2)$, where  $v_2$ is the Veronese embedding; $D$ is the Steiner Roman surface and is homeomorphic to $\mathbb R\mathbb P^2$;  
\subitem-$S= \varphi (\mathbb P^1\times \mathbb P^1)$, where $\varphi$ is the embedding defined by $|l_1+2l_2|$ for $l_{1,2}$ the rulings of $\mathbb P^1\times \mathbb P^1$;
\subitem-$S= \varphi(\mathbb F_2)$, where $\varphi$ is the embedding defined by $|\sigma+ f|$ for $\sigma$ the negative section and $f$ the fibre of $\mathbb F_2$. 
\item[(b)] the image of a surface $\hat D\subset \mathbb P^4$ with canonical singularities under the projection from a point not lying on it;
$\hat D$ is a degenerate dP$4$ surface which is the blowup of $\mathbb P^2$ in $5$ points in almost general position.
\item[(c)]  a rational surface embedded by a complete linear system on its normalisation $\hat D$; the non-normal locus of $D$ is a line $L$ and $D$ may have isolated singularities outside $L$. The minimal resolution of the normalisation of $D$ is a blowup of $\mathbb P^2$ in $9$ points.  
The normalisation of $D$ has at most two rational triple points lying on the inverse image of the non-normal locus; their images on $D$ are also triple points. 
\end{enumerate}
\end{thm}
\begin{rem}
$D$ is not slc in case 1.\end{rem}
\begin{cor}
Let $(X, D_X)$ be a (t,lc) quartic CY pair with non-normal boundary. Then, $(X, D_X)$ has maximal intersection except in the cases described in 2.(a) and (b) of Theorem~\ref{nonnormal4ics}. \end{cor}
\begin{exa} Consider the pair $(X, D_X)$ where: 
\[ X= \{ x_0^2x_3^2+x_1^2x_2x_3+ x_2^2q(x_0,x_1)+ x_4f_3=0\}, D_X= X\cap \{ x_4=0\},\]
where $q$ is a general quadratic form in $(x_1, x_2)$ and $f_3$ a general cubic in $x_0, \cdots, x_4$. 

When $q$ and $f_3$ are general, the quartic hypersurface $X$ has $3$ ordinary double points. Indeed, denote by $L=\{ x_0=x_1=x_4=0\}$, then 
$\Sing (X)$ consists of points of intersection of $L$ with $\{ f_3=0\}$; there are $3$ such points $\{ q_1,q_2, q_3\}$ when $f_3$ is general. 
In the neighbourhood of each point $q_i$ for $i=1,2,3$, the equation of $X$ is of the form 
\[ \{0\}\in \{ xy+zt=0\}\]
(and $D_X= \{t=0\}$) so that all singular points of $X$ are ordinary double points. 
The nodal quartic $X$ is terminal and $\mathbb Q$-factorial because it has less than $9$ ordinary double points; $X$ is birationally rigid by \cite{Che, Me04}. 

Taking the divisorial extraction of the line $L$ is enough to produce a dlt modification $(\widetilde X, D_{\tilde X}+E)$ of $(X, D_X)$; this shows that $(X, D_X)$ does not have maximal intersection. The dual complex has a single 1-stratum, the elliptic curve $D_{\tilde X}\cap E$, which is a $(2,2)$ curve in $\mathbb P^1\times \mathbb P^1$.  The quartic surface $D_X$ is a ruled surface over an elliptic curve isomorphic to $D_{\tilde X}\cap E$; it is an example of case 2.(b) in Theorem~\ref{nonnormal4ics}.
  \end{exa}

\begin{thm}\label{isosingquartic}\cite{Wallp} A normal quartic surface $D\subset \mathbb P^3$ with at least one non-canonical singular point is one of:
\begin{enumerate}[1.]
\item[1.] $D$ has a single elliptic singularity and $D$ is rational, or
\item[2.] $D$ is a cone, or
\item[3.] $D$ is elliptically ruled and 
\item[(a)] $D$ has a double point $p$ with tangent cone $z^2$, the projection away from $p$ is the double cover of $\mathbb P^2$ branched over a sextic curve $\Gamma$. The curve $\Gamma$ is the union of $3$ conics in a pencil that also contains a double line. When this line is a common chord, $D$ has two $T_{2,3,6}$ singularities, when this line is a common tangent, $D$ has one singularity of type $E_{4,0}$. In the first case, $D$ may have an additional $A_1$ singular point. 
\item[(b)] $D$ is $\{ (x_0x_3+ q(x_1, x_2))^2+ f_4(x_1,x_2,x_3)=0\}$ and $\{f_4=0\}$ is four concurrent lines. Depending on whether $L=\{ x_3=0\}$ is one of these lines or not and on whether the point of concurrence lies on $L$, $D$ has either two $T_{2,4,4}$ singular points or one trimodal elliptic singularity. The surface may have additional canonical points $A_n$ for $n=1,2,3$ or $2A_1$.
\end{enumerate}
\end{thm}

\begin{exa}
Let $X$ be the nonsingular quartic hypersurface
\[ X= \{ x_0^2x_3^2+x_0x_1^3+x_3x_2^3+ x_0x_1x_2x_3+ x_4(x_0^3+x_3^3+x_4^3)=0\}\]and $D_X$ its hyperplane section $X\cap \{ x_4=0\}$.
The surface $D_X$ is normal, $$\Sing(D_X)= \{ p,p'\}= \{(1{:}0{:}0{:}0{:}0), (0{:}0{:}0{:}1{:}0)\},$$ and each singular point is simple elliptic $J_{2,0}=T_{2,3,6}$, i.e.~is locally analytically equivalent to $\{0\}\in \{ x^2+y^3+z^6 + xyz=0\}$. 

Here $X$ is nonsingular and $D_X$ is irreducible and normal, and as I explain below, finding a good (t,dlt) modification amounts to constructing a minimal resolution of $D_X$. Let $\widetilde X\to X$ be the composition of the weighted blowups at $p=(1{:}0{:}0{:}0{:}0)$ with weights $(0,2,1,3,1)$ and at $p'=(0{:}0{:}0{:}1{:}0)$ with weights $(3,1,2,0,1)$, and denote by $E$ and $E'$ the corresponding exceptional divisors. Note that $\widetilde X$ is terminal and $\mathbb Q$-factorial by \cite[Theorem 3.5]{Kaw} and has no worse than cyclic quotient singularities. The morphism \[(\widetilde X, D+E+E') \stackrel{f}\to (X, D)\] is volume preserving and the intersection of $D$ with each exceptional divisor is a smooth elliptic curve $C_6\subset \mathbb P(1,1,2,3)$ not passing through the singular points of $E$ and $E'$; $f$ is a good (t,dlt) modification.

The dual complex $\mathcal D(X, D_X)$ is $1$-dimensional, it has $3$ vertices and $2$ edges; $(X, D_X)$ does not have maximal intersection. 
The quartic surface $D_X$ is an example of case 3.(a) in Theorem~\ref{isosingquartic}.
\end{exa}

\begin{cor}
Let $(X,D_X)$ be a (t,lc) quartic CY pair. Assume that $D_X$ is normal, has non-canonical singularities but is not a cone. Then $(X,D_X)$ has maximal intersection except in cases  3.(a) and (b) of Theorem~\ref{isosingquartic}.
\end{cor}
\begin{rem} When $\dim \mathcal D (X, D_X)=1$, $D_X$ either has two $T_{2,3,6}$ or two $T_{2,4,4}$ singularities. Indeed, as is explained in Section~\ref{sings}, singular points $p\in D$ are Kodaira singularities, and in particular are at worst bimodal. The description of cases 3.(a) and (b) of Theorem~\ref{isosingquartic} immediately implies the result, because a surface singularity of type $E_{4,0}$ is trimodal. \end{rem}

\end{document}